\documentclass[11pt]{article}
\usepackage{amssymb}
\usepackage{amsmath}
\usepackage{graphicx}

\title{ }
\author{}
\date{}

\newtheorem{theorem}{Theorem}[section]

\newtheorem{coro}[theorem]{Corollary}

\newtheorem{prob}[theorem]{Problem}

\newcommand{\qed}{\hspace*{\fill} \rule{7pt}{7pt}}

\begin{document}


\title{On edges not in monochromatic copies of a fixed bipartite graph}

\author{Jie Ma\thanks{School of Mathematical Sciences, University of Science and Technology of China,
Hefei, Anhui 230026, China. Email: jiema@ustc.edu.cn. Partially supported by NSFC project 11501539.}}

\maketitle

\begin{abstract}
Let $H$ be a fixed graph. Denote $f(n,H)$ to be the maximum number of edges not contained in any monochromatic
copy of $H$ in a 2-edge-coloring of the complete graph $K_n$, and $ex(n,H)$ to be the {\it Tur\'an number} of $H$.
An easy lower bound shows $f(n,H)\ge ex(n,H)$ for any $H$ and $n$.
In \cite{KS2}, Keevash and Sudakov proved that if $H$ is an edge-color-critical graph or $C_4$,
then $f(n,H)= ex(n,H)$ holds for large $n$,
and they asked if this equality holds for any graph $H$ when $n$ is sufficiently large.
In this paper, we provide an affirmative answer to this problem for an abundant
infinite family of bipartite graphs $H$, including all even cycles and
complete bipartite graphs $K_{s,t}$ for $t>s^2-3s+3$ or $(s,t)\in\{(3,3),(4,7)\}$.
In addition, our proof shows that for all such $H$, the 2-edge-coloring $c$ of $K_n$ achieves
the maximum number $f(n,H)$ if and only if one of the color classes in $c$ induces an extremal graph
for $ex(n,H)$. We also obtain a multi-coloring generalization for bipartite graphs.
Some related problems are discussed in the final section.
\end{abstract}

\section{Introduction}
Given a graph $H$, let $f(n,H)$ be the maximum number of edges not contained in any monochromatic
copy of $H$ in a 2-edge-coloring of the complete graph $K_n$, and let $ex(n,H)$ be the {\it Tur\'an number} of $H$,
i.e., the maximum number of edges in an $n$-vertex $H$-free graph.
The problem of determining $f(n,H)$ was motivated by counting the number of monochromatic cliques,
and we refer interested readers to \cite{KS2} for a thoughtful discussion on the background and related topics.
(For results on monochromatic cliques, see \cite{G59,T97,E97,EFGJL,KS1,Y08,CKPSTY}.)

If one considers the 2-edge-coloring of $K_n$ in which one of the colors induces the largest $H$-free graph,
then it is easy to see that for any $H$ and $n$, we have
\begin{align}\label{equ:f>=ex}
f(n,H)\ge ex(n,H).
\end{align}
Erd\H{o}s, Rousseau and Schelp (see \cite{E97}) showed that $f(n,K_3)=ex(n,K_3)$ for sufficiently large $n$,
and this also can be derived from a result of Pyber in \cite{P86} for $n\ge 2^{1500}$.
The generalization of this result was suggested by Erd\H{o}s in \cite{E97}.
Keevash and Sudakov \cite{KS2} studied general graphs and asked that if,
for large $n$, the above lower bound \eqref{equ:f>=ex} is tight.
\begin{prob} (\cite{KS2})\label{prob:main}
Let $H$ be a fixed graph. Is it true that for $n$ sufficiently large, $f(n,H)=ex(n,H)$?
\end{prob}
The authors of \cite{KS2} confirmed it for $H$ being any edge-color-critical graph or a $C_4$,
and in fact, quite amazingly, they were able to determine the value of $f(n,H)$ for every $n$ when $H$ is a $K_3$ or $C_4$.
We quote from their remark \cite{KS2} that ``for bipartite graphs the situation is less clear,
as even the asymptotics of the Tur\'an numbers are known only in a few cases".

In this paper, we provide an affirmative answer to Problem \ref{prob:main} for an abundant infinite family of bipartite graphs.
A vertex $w$ in a bipartite graph $H$ is called {\it weak}, if
$$ex(n,H-w)=o(ex(n,H)).$$
The notation of weak vertices is explicitly defined in the literature and has been well studied (see \cite{S84}).
We call a bipartite graph $H$ {\it reducible}, if it contains a weak vertex $w$
such that $H-w$ is connected. For instance all even cycles are reducible.
Our main theorem is as follows.
\begin{theorem}\label{thm:main}
Let $H$ be a reducible bipartite graph. Then for sufficiently large $n$, $f(n,H)=ex(n,H)$.
Moreover, a 2-edge-colorings of $K_n$ achieves the maximum number $f(n,H)$
if and only if one of the color classes induces an extremal graph for $ex(n,H)$.
\end{theorem}
We point out that the ``moreover" part is new for $C_4$,
while its analog is not true for edge-color-critical graphs as noticed in \cite{KS2}.

Let $\mathcal{C}^*$ be the family of bipartite graphs,
each of which contains a cycle and a vertex $w$ whose deletion will result in a tree.
It is easy to see that all graphs in $\mathcal{C}^*$, including even cycles and Theta graphs,\footnote{The Theta
graphs $\theta_{k,l}$ denotes the graph consisting of $k$ internally disjoint paths of length $l$ between
two fixed endpoints, for $k,l\ge 2$.} are reducible.
Based on the current knowledge on degenerated Tur\'an numbers,
we collect some reducible graphs in the coming result.
\begin{coro}\label{Coro}
For $n$ sufficiently large, $f(n,H)=ex(n,H)$ holds for every $H$ as following: even cycles $C_{2l}$,
Theta graphs $\theta_{k,l}$,
and complete bipartite graphs $K_{s,t}$ for $t>s^2-3s+3$ or $(s,t)\in\{(3,3),(4,7)\}$.
\end{coro}

The rest of the paper is organized as follows. In the next section, we prove Theorem \ref{thm:main}
in full and then derive Corollary \ref{Coro}. In Section 3, we generalize Theorem \ref{thm:main} to multi-colorings.
In the final section, we close this paper by mentioning some related problems.

\section{Reducible bipartite graphs}
Let $H$ be a fixed graph and $c$ be a $k$-edge-coloring of $K_n$.
An edge of $K_n$ is called {\it NIM-H},
if it is not contained in any monochromatic copy of $H$ in $c$.
Let $E_c$ denote the set of all NIM-H edges of $K_n$.
For $A,B\subseteq V(K_n)$, by $(A,B)$ we denote the complete bipartite graph with two parts $A$ and $B$.

In this section, we establish Theorem \ref{thm:main} and Corollary \ref{Coro}.
To do so, we prove the following stronger result.
\begin{theorem}\label{thm:str}
Let $H$ be a reducible bipartite graph.
If $c$ is a 2-edge-coloring of $K_n$ such that $E_c$ contains a red edge and a blue edge,
then $|E_c|=o(ex(n,H))$.
\end{theorem}
\begin{pf}
Let $h=|V(H)|$, $(X,Y)$ be the bipartition of $H$, and $w\in X$ be a weak vertex of $H$
such that $ex(n,H-w)=o(ex(n,H))$ and $H-w$ is connected. Note that $(X-w,Y)$ is the unique bipartition of $H-w$, as $H-w$ is connected.

We first define a red star $S_0$ in $K_n$ (i.e., all edges in the star are red), which contains at least one red NIM-H edge, as follows.
If there exist vertices incident with a red NIM-H edge and at least $h$ red edges,
then pick one such vertex $x$ and form a star $S_0$ consisting of the center $x$ and $h$ red neighbors of $x$
such that $xv$ is a red NIM-H edge for some $v\in V(S_0)$.
Otherwise every vertex incident with a red NIM-H edge has less than $h$ red neighbors,
then pick one such vertex $x$ with maximum number of red neighbors and let $S_0$ consist of $x$ and all its red neighbors.

Similarly as above, we define a blue star $S_1$ in $K_n$, which contains at least one blue NIM-H edge.
Let $x, y$ be the centers of the stars $S_0, S_1$, respectively.
Note that $S_0$ and $S_1$ may share some common vertices. We let $$S=V(S_0)\cup V(S_1)=\{s_1,s_2,...,s_t\}.$$
So $t=|S|\le 2h+2$. For $z\in V(K_n)\backslash S$, let $\vec{\epsilon}(z)=(\epsilon_1,\epsilon_2,...,\epsilon_t)$ be the vector
such that
\begin{equation}
\epsilon_i=\left\{\begin{array}{ll}
    0 & \text{if } zs_i \text{ is red}, \\
    1 & \text{if } zs_i \text{ is blue}.
  \end{array}\right.
\end{equation}
For $\vec{v}\in \{0,1\}^t$, let $A_{\vec{v}}$ denote the set of all vertices $z\in V(K_n)\backslash S$
such that $\vec{\epsilon}(z)=\vec{v}$.
Observe that all edges between $s_i\in S$ and $A_{\vec{v}}$ must be monochromatic.

We now consider the numbers of NIM-H edges adjacent to sets $A_{\vec{v}}$.
The first claim implies that the number of NIM-H edges adjacent to $A_{\vec{0}}\cup A_{\vec{1}}$ is $O(n)$.

\medskip

{\bf Claim 1}: $|A_{\vec{0}}|<h$ and $|A_{\vec{1}}|<h$.

By symmetry, it suffices to consider $A_{\vec{0}}$.
We notice that all edges in $(A_{\vec{0}},S)$ are red.
Suppose for a contradiction that $|A_{\vec{0}}|\ge h$.
If the red star $S_0$ has less than $h+1$ vertices, then it is clear
that no vertex in $V(K_n)\backslash S$ can be adjacent to $x$, implying that $A_{\vec{0}}=\emptyset$.
So the red star $S_0$ has exactly $h+1$ vertices.
We see that all edges in $(A_{\vec{0}}\cup \{x\},S_0-\{x\})$ are red.
From this, one can easily find a red copy of $H$ which uses one NIM-H edge of $x$,
contradicting the definition of NIM-H edges. This proves claim 1.

\medskip

{\bf Claim 2}: For $\vec{v}\in \{0,1\}^t-\{\vec{0},\vec{1}\}$, the number of NIM-H edges contained in $A_{\vec{v}}$ is at most $2\cdot ex(n,H-w)$.

As $\vec{v}\notin \{\vec{0},\vec{1}\}$, there exist $a,b\in S$ such that all edges
in $(a,A_{\vec{v}})$ are red and all edges in $(b,A_{\vec{v}})$ are blue.
If the red NIM-H edges in $A_{\vec{v}}$ form a copy $K=(X-w,Y)$ of $H-w$,
then $K\cup \{a\}$ would contain a red copy of $H$ with some NIM-H edges, a contradiction.
Therefore, neither the red NIM-H edges nor the blue NIM-H edges can form a copy of $H-w$.
This proves claim 2.

\medskip

{\bf Claim 3}: For $\vec{v},\vec{u}\in \{0,1\}^t-\{\vec{0},\vec{1}\}$, the number of NIM-H edges in
$(A_{\vec{v}},A_{\vec{u}})$ is at most $2\cdot ex(n,H-w)$.

Suppose that the red NIM-H edges in $(A_{\vec{v}},A_{\vec{u}})$ form a copy $K=(X-w,Y)$ of $H-w$.
By symmetry, we assume that $X-w\subseteq A_{\vec{v}}$ and $Y\subseteq A_{\vec{u}}$.
Since $\vec{u}\neq \vec{1}$, there exists $a\in S$ such that all edges in $(a,A_{\vec{u}})$ are red.
Adding $a$ and all red edges in $(a,Y)$ to $K$ would result in a red copy of $H$, a contradiction.
Therefore, neither the red NIM-H edges nor the blue NIM-H edges in $(A_{\vec{v}},A_{\vec{u}})$ can form a copy of $H-w$.
Claim 3 is finished.

\medskip

Each edge in $E_c$ is either adjacent to $S\cup A_{\vec{0}}\cup A_{\vec{1}}$
or contained in $A_{\vec{v}}$ or $(A_{\vec{v}},A_{\vec{u}})$ for some $\vec{v},\vec{u}\in \{0,1\}^t-\{\vec{0},\vec{1}\}$.
Since $|S|=t\le 2h+2$, there are at most $2^{2h+2}$ sets $A_{\vec{v}}$. Combining the above claims, we have
\begin{align*}
|E_c|&\le (|S|+2h)\cdot n +2^{2h+2}\cdot 2\cdot ex(n,H-w)+(2^{2h+2})^2\cdot 2\cdot ex(n,H-w)\\
&\le 2^{4h+6}\cdot ex(n,H-w)=o(ex(n,H)).
\end{align*}
This finishes the proof of Theorem \ref{thm:str}.
\qed
\end{pf}

\bigskip

We are ready to prove Theorem \ref{thm:main}.

\medskip

\noindent{\it Proof of Theorem \ref{thm:main}.}
We have seen $f(n,H)\ge ex(n,H)$ from \eqref{equ:f>=ex}.

Let $c$ be a 2-edge-coloring of $K_n$ such that $|E_c|=f(n,H)$.
If $E_c$ contains a red edge and a blue edge, then by Theorem \ref{thm:str},
we have $ex(n,H)\le f(n,H)=|E_c|=o(ex(n,H))$, a contradiction.
So we may assume that all NIM-H edges are red.
It then becomes clear that $E_c$ does not contain any copy of $H$,
implying that $f(n,H)=|E_c|\le ex(n,H)$. This proves that $f(n,H)=ex(n,H)$ for large $n$.

It also follows that $|E_c|=ex(n,H)$. So $E_c$ must induce an extremal graph for $ex(n,H)$.
We claim that except these edges in $E_c$, no other edge can be red.
Suppose not, say $e\in E(K_n)-E_c$ is red.
Then $E_c\cup \{e\}$ induces an $n$-vertex graph with more than $ex(n,H)$,
which must contain a copy of $H$. But this $H$ contains all red edges
and in particular some NIM-H edges from $E_c$, a contradiction. This proves the claim.
Now we see that all red edges of $c$ induces an extremal graph for $ex(n,H)$.

To prove the ``moreover" part, it remains to show that
if all red edges of $c$ induces an extremal graph for $ex(n,H)$, then $|E_c|=f(n,H)=ex(n,H)$.
Since all red edges surely are NIM-H, we have $|E_c|\ge ex(n,H)$.
So we need to show that no blue edge can be NIM-H. This, again, can be derived from Theorem \ref{thm:str}.
We have finished the proof. \qed

\bigskip

We conclude this section by showing Corollary \ref{Coro}.
Recall the seminal theorem of K\H{o}v\'ari-S\'os-Tur\'an \cite{KST} and
the best known general lower bound on Tur\'an number of $K_{s,t}$ that
$$\Omega(n^{2-\frac{s+t-2}{st-1}})\le ex(n,K_{s,t})\le \frac12(t-1)^{1/s}n^{2-1/s}+\frac12(s-1)n.$$
We also need the result (see \cite{KRS,ARS}) that $ex(n,K_{s,t})\ge \Omega(n^{2-1/s})$ for $t>(s-1)!$.

\medskip

\noindent{\it Proof of Corollary \ref{Coro}.}
In view of Theorem \ref{thm:main}, it is enough to show that every graph $H$ in the list is reducible.
As it is clear that even cycles and Theta graphs are reducible, we only need to consider $K_{s,t}$.
When $t>(s-1)!$, it holds that $$\frac{ex(n,K_{s-1,t})}{ex(n,K_{s,t})}=O\left(\frac{n^{2-1/(s-1)}}{n^{2-1/s}}\right)=o(1),$$
and when $t>s^2-3s+3$, we have $\frac{s+t-2}{st-1}<\frac{1}{s-1}$, implying that
$$\frac{ex(n,K_{s-1,t})}{ex(n,K_{s,t})}= O\left(\frac{n^{2-1/(s-1)}}{n^{2-(s+t-2)/(st-1)}}\right)=o(1).$$
Therefore, $K_{s,t}$ is reducible whenever $t>\min\{s^2-3s+3,(s-1)!\}$, finishing the proof.\qed

\section{Generalization to multi-colorings}
In this section, we consider multi-color versions of Theorem \ref{thm:main}.

For $k\ge 3$, let $f_k(n,H)$ denote the maximum number of edges not contained in any monochromatic
copy of $H$ in a $k$-edge-coloring of the complete graph $K_n$.
Given a bipartition $(X,Y)$ of bipartite $H$, let $ex^*(m,n,H)$ denote the maximum number of edges of
graphs $G$, where $G$ is a spanning subgraph of $K_{m,n}$ and has no copies of $H=(X,Y)$ with $X$ contained in the $m$-part.\footnote{The Zarankiewicz function $z(m,n,s,t)$ is just the same as $ex^*(m,n,K_{s,t})$.}

We first prove a general lower bound for every bipartite graph $H$ that
\begin{align}\label{equ:fk>=ex}
f_k(n,H)\ge (k-1)\cdot ex(n,H)- O\left(ex(n,H)/n\right)^2=(k-1-o(1))\cdot ex(n,H).
\end{align}

\noindent{\it Proof.}
Let $G$ be an $n$-vertex $H$-free extremal graph for $ex(n,H)$.
For a permutation $\pi$ on $V(G)$, let $G(\pi)$ be obtained from $G$
by permuting all edges according to $\pi$, i.e., $E(G(\pi))=\pi(E(G))$.
Take $k-1$ random permutations $\pi_1,\pi_2,...,\pi_{k-1}$ and consider the overlap $E_{ij}=E(G(\pi_i))\cap E(G(\pi_j))$.
Since the probability that each $e\in \binom{V}{2}$ belongs to $G(\pi_i)$ equals $ex(n,H)/\binom{n}{2}$,
the expectation of $\sum_{i,j} |E_{ij}|$ is at most $\binom{k-1}{2} ex(n,H)^2/\binom{n}{2}$.
Therefore, there exist permutations $\pi_1,\pi_2,...,\pi_{k-1}$ such that
the total overlap $\sum_{i,j} |E_{ij}|$ is at most $\binom{k-1}{2} ex(n,H)^2/\binom{n}{2}$.
We then define a $k$-edge-coloring $c$ of $K_n$ as following.
Color the edges of $G(\pi_1)$ by color 1; and for $2\le i\le k-1$,
color the edges in $E(G(\pi_i))-\cup_{1\le j\le i-1} E_{ij}$ by color $i$;
and lastly, color all edges of $K_n$ not in $\cup_{1\le i\le k-1} E(G(\pi_i))$ by color $k$.
This implies that
$f_k(n,H)\ge |\cup_{1\le i\le k-1} E(G(\pi_i))|\ge (k-1)\cdot ex(n,H)-\sum_{i,j}|E_{ij}|,$
which is at least $(k-1)\cdot ex(n,H)- O\left(ex(n,H)/n\right)^2$. \qed

\begin{theorem}\label{thm:k-C4}
For $n$ sufficiently large, we have $$(k-1-o(1))\cdot ex(n,C_4)\le f_k(n,C_4)\le (k-1)\cdot ex(n,C_4).$$
\end{theorem}

\begin{theorem}\label{thm:multi}
Let $H$ be a bipartite graph with a vertex $w$ such that
$ex^*(n,n,H-w)=o(ex(n,H)).$
Then for sufficiently large $n$, $$(2-o(1))\cdot ex(n,H)\le f_3(n,H)\le 2\cdot ex(n,H).$$
Such graphs $H$ include even cycles $C_{2l}$ and
complete bipartite graphs $K_{s,t}$ for $t>s^2-3s+3$ or $(s,t)\in\{(3,3),(4,7)\}$.
\end{theorem}

\begin{pf} (For both Theorems \ref{thm:k-C4} and \ref{thm:multi}.)
The lower bound follows from $\eqref{equ:fk>=ex}$.

First we prove an analog of Theorem \ref{thm:str}.
Let $(X,Y)$ be the partition of $H$ which $ex^*(n,n,H-w)$ refers to. Let $w\in X$ and $h=|V(H)|$.
Call an edge with color $i$ as an $i$-edge for convenience.

{\bf Claim:} Let $c$ be a $k$-edge-coloring of $K_n$. If $E_c$ contains a NIM-$H$ $i$-edge
for each $i\in [k]$, then $|E_c|\le (k-2)\cdot ex(n,H)+o(ex(n,H))$.

The proof of this claim will follow the same lines of Theorem \ref{thm:str}.
For each color $i\in [k]$, we define a star $S_i$ in $K_n$ consisting of $i$-edges,
among which there is at least one NIM-H $i$-edge.
If there exist vertices incident with a NIM-H $i$-edge and at least $h$ $i$-edges,
then pick one such vertex $x_i$ and form a star $S_i$ with the center $x_i$ and consisting of $h$ $i$-edges
such that there exists at least one NIM-H $i$-edge $x_iv$ for some $v\in V(S_i)$.
Otherwise every vertex incident with a NIM-H $i$-edge has less than $h$ $i$-neighbors,
then pick one such vertex $x_i$ with maximum number of $i$-neighbors and let $S_i$ consist of $x_i$ and all its $i$-neighbors.
Let $S=\cup_{i\in [k]}V(S_i)=\{s_1,s_2,...,s_t\}.$
So $t=|S|\le k(h+1)$. For $z\in V(K_n)\backslash S$,
let $\vec{\epsilon}(z)=(\epsilon_1,\epsilon_2,...,\epsilon_t)$ be the vector
such that $\epsilon_i=j$ iff $zs_i$ is colored by $j$.
For $\vec{v}\in [k]^t$, let $A_{\vec{v}}$ denote the set of all vertices $z\in V(K_n)\backslash S$
such that $\vec{\epsilon}(z)=\vec{v}$.

For some $I\subseteq [k]$, we say $A_{\vec{v}}$ is {\it $I$-feasible}, if for each $i\in I$ there exists some coordinate in $\vec{v}$ being $i$,
and subject to this, $I$ is maximal. We then establish the following three assertions.

\medskip

{\bf (1).} For each $i\in [k]$, we have $|A_{\vec{i}}|<h$.

Note that all edges in $(A_{\vec{i}}, S)$ are of color $i$. If $|A_{\vec{i}}|\ge h$,
then the complete bipartite graph $(A_{\vec{i}}\cup \{x_i\}, S_i-\{x_i\})$ contains
a copy of $H$ of all $i$-edges with at least one NIM-H $i$-edge (incident to $x_i$),
a contradiction. This proves (1).

\medskip

{\bf (2).} For $i\in I$, the $I$-feasible set $A_{\vec{v}}$ has no more than $ex(n,H-w)$ NIM-H $i$-edges.

Suppose that the NIM-H $i$-edges in $A_{\vec{v}}$ form a copy $K$ of $H-w$.
Since $i\in I$, there exists some $a\in S$ such that the edges
in $(a, A_{\vec{v}})$ are all of color $i$.
Then adding $a$ into $K$ would give a copy of $H$ of color $i$
which also contains NIM-H edges, a contradiction.
This shows that there are no more than $ex(n,H-w)$ NIM-H $i$-edges in $A_{\vec{v}}$, establishing (2).

\medskip

{\bf (3).} Let $A_{\vec{u}}$ be $I$-feasible and $A_{\vec{v}}$ be $J$-feasible.
For $i\in I\cup J$, there are no more than $ex^*(n,n,H-w)$ NIM-H $i$-edges in $(A_{\vec{u}},A_{\vec{v}})$.

For $i\in I\cup J$, there exists some coordinate in $\vec{u}$ or $\vec{v}$ being $i$.
By symmetry, say this coordinate is from $\vec{v}$. Then there exists some $a\in S$
such that all edges from $a$ to $A_{\vec{v}}$ are of color $i$.
Suppose that the NIM-H $i$-edges in $(A_{\vec{u}},A_{\vec{v}})$
contains a copy $K$ of $H-w=(X-w,Y)$ with $X-w\subseteq A_{\vec{u}}$ and $Y\subseteq A_{\vec{v}}$.
Then $\{a\}\cup K$ would contain a copy of $H$ of color $i$ with some NIM-H edges,
a contradiction. Thus $(A_{\vec{u}},A_{\vec{v}})$ has no more than
$ex^*(|A_{\vec{u}}|,|A_{\vec{v}}|,H-w)\le ex^*(n,n,H-w)$ NIM-H $i$-edges. This proves (3).

\medskip

Observe that the NIM-H edges not in (2) and (3)
are of the following three types:
\begin{itemize}
\item[(i).] NIM-H edges which are adjacent to $S$ or $A_{\vec{i}}$ for some $i\in [k]$,
\item[(ii).] NIM-H $i$-edges in $I$-feasible set $A_{\vec{v}}$, where $i\in I^c$ and $|I|\ge 2$, and
\item[(iii).] NIM-H $i$-edges between $I$-feasible set $A_{\vec{u}}$ and $J$-feasible set $A_{\vec{v}}$,
where $i\notin I\cup J$ (or equivalently $i\in I^c\cap J^c$) and $|I|,|J|\ge 2$.
\end{itemize}
Note that $|S|=t\le k(h+1)$. So there are at most $k^{k(h+1)}$ sets $A_{\vec{v}}$, which is constantly many.
Also note that $ex(n,H-w)=O(ex^*(n,n,H-w))=o(ex(n,H))$,
implying that $H$ must not be a forest and thus $ex(n,H)=\Omega(n^{1+c})$ for some $c>0$.
These, combining with the above assertions, imply that
there are just $o(ex(n,H))$ NIM-H edges contained in (2), (3) and (i).
To complete the proof of the claim, it then suffices to show that the number $N^*$ of NIM-H edges
in (ii) and (iii) is at most $(k-2)\cdot ex(n,H)+o(ex(n,H))$.

For $i\in [k]$, denote $B_i$ to be the union of all $I$-feasible sets $A_{\vec{v}}$ satisfying $i\in I^c$ and $|I|\ge 2$.
It is straightforward to verify that the NIM-H $i$-edges in (ii) and (iii) must be contained in $B_i$.
Thus the number of NIM-H edges in (ii) and (iii) is $N^*\le \sum_{i=1}^k ex(b_i,H),$ where $b_i=|B_i|$.
Note each $b_i\le n$ and $\sum_{i=1}^k b_i\le (k-2)\cdot\sum |A_{\vec{v}}|\le (k-2)n$,
as every vertex in $I$-feasible sets with $|I|\ge 2$ only can appear in at most $k-2$ $B_i$'s.

When $k=3$, we have $b_1+b_2+b_3\le n$, so it is clear that $N^*\le \sum_{i=1}^3 ex(b_i,H)\le ex(n,H).$
When $H=C_4$, using the well-known result (see \cite{ERS,B66,KST}) that $ex(n,C_4)=(1/2+o(1))\cdot n^{3/2}$,
it holds that $N^*\le \sum_{i=1}^k ex(b_i,C_4)\le \frac{1}{2}\sum_{i=1}^k b_i^{3/2}+o(n^{3/2})$.
Subject to $b_i\le n$ and $\sum_{i=1}^k b_i\le (k-2)n$, by convexity, we have
$N^*\le \frac{k-2}{2}\cdot n^{3/2}+o(n^{3/2})=(k-2)\cdot ex(n,C_4)+o(ex(n,C_4))$, which is desired.
This completes the proof of the claim.

\medskip

Next we prove the upper bound of $f_k(n,H)$.
Let $c$ be a $k$-edge-coloring of $K_n$ such that $|E_c|=f_k(n,H)$.
If $E_c$ contains a NIM-H $i$-edge for every $i\in [k]$,
then by the claim, we have
$(k-1-o(1))\cdot ex(n,H)\le f_k(n,H)=|E_c|\le (k-2)\cdot ex(n,H)+o(ex(n,H))$, a contradiction.
Therefore $E_c$ has at most $k-1$ colors.
The set of NIM-H edges of the same color contains none copy of $H$ and thus is of size at most $ex(n,H)$.
Thus, $f_k(n,H)=|E_c|\le (k-1)\cdot ex(n,H)$.

It remains to verify that $ex^*(n,n,H-w)=o(ex(n,H))$ holds for $H$ being an even cycle or $K_{s,t}$
for $t>\min\{s^2-3s+3,(s-1)!\}$. This follows by the same proof of Corollary \ref{Coro},
using $ex^*(m,n,K_{s,t})\le (t-1)^{1/s}mn^{1-1/s}+(s-1)n$ (see \cite{KST}).
This proves Theorems \ref{thm:k-C4} and \ref{thm:multi}. \qed
\end{pf}

\section{Concluding remarks}
In Theorem \ref{thm:main} we prove that $f(n,H)$ equals $ex(n,H)$ for bipartite graphs $H$
having weak vertices (for sufficiently large $n$).
Simonovits asked in \cite{S84} to ``characterize those bipartite graphs which have weak vertices" and this remains unclaimed.

It seems that for $k\ge 3$, the function $f_k(n,H)$ has a different behavior between bipartite and non-bipartite graphs.
For bipartite $H$, it may be reasonable to ask if $f_k(n,H)=(k-1)\cdot ex(n,H)$ holds for sufficiently large $n$.
For non-bipartite graphs, the situation is more complicate.
We speculate that the following 3-edge-coloring $c$ of $K_n$ (which also is the extremal configuration as in \cite{CKPSTY})
achieves the maximum of $f_3(n,K_3)$:
let $V(K_n)=V_1\cup V_2\cup...\cup V_5$, where $||V_i|-|V_j||\le 1$, and
color all edges in $(V_i,V_{i+1})$ by red, all edges in $(V_i,V_{i+2})$ by blue and all edges in each $V_i$ by green.
For more discussion and other related problems, we direct readers to \cite{KS2}.



\end{document}